\documentclass{ps}

\usepackage[german,english]{babel}
\usepackage[dvips]{graphics}

\usepackage[T1]{fontenc}
\usepackage[latin1]{inputenc}

\def \beq {\begin{eqnarray}}
\def \eeq {\end{eqnarray}}
\def \beqn {\begin{eqnarray*}}
\def \eeqn {\end{eqnarray*}}

\newcounter{for}[section]
\newtheorem{itremark}{Remark}[section]
\newenvironment{remark}{\begin{itremark}\rm}{\end{itremark}}
\newcommand{\be}[1]{\addtocounter{for}{1} \begin{equation}\label{#1}}
\newcommand{\ee}{\end{equation}}
\newcommand{\br}[1]{\begin{remark}\label{#1}}

\newcommand{\bi}{\begin{itemize}}
\newcommand{\ei}{\end{itemize}}
\newcommand{\ben}{\begin{enumerate}}
\newcommand{\een}{\end{enumerate}}

\def\vbar{\mathchoice{\vrule height6.3ptdepth-.5ptwidth.8pt\kern-.8pt}
   {\vrule height6.3ptdepth-.5ptwidth.8pt\kern-.8pt}
   {\vrule height4.1ptdepth-.35ptwidth.6pt\kern-.6pt}
   {\vrule height3.1ptdepth-.25ptwidth.5pt\kern-.5pt}}
\def\fudge{\mathchoice{}{}{\mkern.5mu}{\mkern.8mu}}
\def\bbc#1#2{{\rm \mkern#2mu\vbar\mkern-#2mu#1}}
\def\bbb#1{{\rm I\mkern-3.5mu #1}}
\def\bba#1#2{{\rm #1\mkern-#2mu\fudge #1}}
\def\bb#1{{\count4=`#1 \advance\count4by-64 \ifcase\count4\or\bba A{11.5}\or
   \bbb B\or\bbc C{5}\or\bbb D\or\bbb E\or\bbb F \or\bbc G{5}\or\bbb H\or
   \bbb I\or\bbc J{3}\or\bbb K\or\bbb L \or\bbb M\or\bbb N\or\bbc O{5} \or
   \bbb P\or\bbc Q{5}\or\bbb R\or\bbc S{4.2}\or\bba T{10.5}\or\bbc U{5}\or
   \bba V{12}\or\bba W{16.5}\or\bba X{11}\or\bba Y{11.7}\or\bba Z{7.5}\fi}}

\def \Z {{\bb Z}}
\def \R {{\bb R}}

\def \ra {\rightarrow }

\def \s {\sigma}

\def \e {\epsilon}

\def \G {{\cal{G}}}

\def \L {\Lambda}
\def \l {\lambda}
\def \b {\beta}
\def \g {\gamma}

\begin{document}


\title{Stationary measures and phase transition for a class of
Probabilistic Cellular Automata}
\author{Paolo Dai Pra}
\address{Dipartimento di Matematica Pura e %
  Applicata, Universit\`a di Padova, Via Belzoni 7, 35131 Padova, Italy %
\email{daipra@math.unipd.it}}
\author{Pierre-Yves Louis}
\address{Laboratoire de Statistique et Probabilit\'es, %
F.R.E. C.N.R.S. 2222, U.F.R. de Math\'{e}matiques, %
Universit\'{e} Lille 1, 59655 Villeneuve d'Ascq Cedex, France, %
\email{louis@lps.univ-lille1.fr}}
\author{Sylvie R\oe lly}
\address{\foreignlanguage{german}{Weierstra"s-Institut f\"ur Angewandte Analysis
und Stochastik}, %
 \foreignlanguage{german}{Mohrenstra"se 39}, %
 10117 Berlin, Germany, %
\email{roelly@wias-berlin.de}}
\secondaddress{Permanent %
  address: CMAP, UMR CNRS 7641 \'Ecole Polytechnique, 91128 Palaiseau %
  Cedex France.}

\date{Received December 4, 2001. Revised March 11, 2002.}

\begin{abstract}
We discuss various properties of Probabilistic Cellular Automata, such
as the structure of the set of stationary measures and multiplicity of
stationary measures (or phase transition) for reversible models.
\end{abstract}



\maketitle



\section{Introduction}

Probabilistic Cellular Automata (PCA) are discrete-time Markov
chains on a product space $S^{\L}$ ({\em configuration space})
whose transition probability is a product measure. In this paper,
$S$ is assumed to be a finite set ({\em spin space}), and $\L$
(set of {\em sites}) a subset, finite or infinite, of $\Z^{d}$.
The fact that the transition probability $P(d\s|\s')$, $\s,\s' \in
S^{\L}$, is a product measure means that all spins $\{\s_{i}: i
\in \L\}$ are simultaneously and independently updated ({\em
parallel updating}). This transition mechanism differs from the
one in the most common Gibbs samplers (e.g. \cite{gu,br}), where only
one site is updated at each time step ({\em sequential updating}).

Several properties of PCA's, mainly of general and qualitative
nature, have been investigated (\cite{lms,gklm,tvs,daw,ms2}). As
far as we know, however, sharper properties like e.g. rate of
convergence to equilibrium or use of parallel dynamics in perfect sampling,
have not yet been investigated. PCA's are hard to analyze mainly
for the following reason. Suppose $\L$ is a finite subset of
$\Z^{d}$, and let $\mu$ be a given probability on $S^{\L}$. To fix
ideas, we may think of $\mu$ as a finite volume Gibbs measure for
a given interaction and assigned boundary conditions. It is simple
to construct Markov chains on $S^{\L}$ with sequential updating which
have $\mu$ as reversible measure. Transition probabilities are
given in simple form in terms of $\mu$, and reversibility
immediately implies that $\mu$ is an stationary measure for the
dynamics. Quite differently, for a given $\mu$, there
is no general recipe to construct a PCA for which $\mu$ is
stationary. In particular, there exists Gibbs measures on $S^{\Z^{2}}$
such that no PCA admits them  as stationary measures (Theorem 4.2 in
\cite{daw}).

Despite of this descouraging starting point, other aspects of PCA's
make them interesting stochastic models, and motivate further
investigation.
\ben
\item
For simulation and sampling, PCA's are natural stochastic algorithms
for parallel computing. At least in some simple models (see Section
3) it is interesting to evaluate their performance versus
algorithms with sequential updating. This will be the subject of a
forthcoming paper.
\item
In opposition to dynamics with sequential updating, it is simple to
define PCA's in infinite volume without passing to
continuous time. One may try to study, for instance, convergence
to equlibrium in infinite volume, or in finite volume uniformly in the volume size.
Although some perturbative methods are available (see \cite{mami}
Chapter 7, \cite{ms1,ms2}), a theory corresponding the one in (\cite{ma}) in
continuous time, is yet to be developed.
\item
PCA's that are reversible with respect to a Gibbs measure $\mu$
have been completely characterized in \cite{kv}. In particular it
has been shown that only a small class of Gibbs measures may be
reversible for a PCA. For such PCA's one can investigate metastable
behavior.
A first step in this direction is done in \cite{bcls}. \een

The present paper is a small step toward a better understanding of PCA's. Our
objective is first to present some links between the sets of reversible,
resp. stationary, resp. Gibbs measures for general PCA's.
 We then illustrate these results on a particular class of reversible
PCA's already introduced in \cite{bcls}.

More precisely it was proved in \cite{kv} that for PCA's possessing
a reversible
Gibbs measure w.r.t. a potential $\Phi$, all reversible measures are gibbsian
w.r.t the same potential.
We prove a similar statement on the set of stationary measures :
For a general PCA, if one shift invariant
stationary measure is
Gibbsian for a potential $\Phi$, then all shift invariant
stationary measures are Gibbsian
 w.r.t. the same potential $\Phi$ (see Proposition \ref{p2.2}).
This induce that for a class of local, shift invariant, non-degenerated,
reversible PCA the reversible measures coincide  with the  Gibbsian
stationary ones (Remark \ref{r3.2}).

Applying this general statements to the class of PCA's considered in
 \cite{bcls}, one can do explicit a stationary measure
which is in fact Gibbsian w.r.t. a certain potential $\Phi$ we write
 down (cf Proposition \ref{p3.2}); we show that, for sufficiently
small values of the temperature parameter, {\em phase transition} occurs,
that is there are several Gibbs
measures w.r.t. $\Phi$. At least in
certain cases, existence of phase transition would follow from
general expansion arguments, like Pirogov-Sinai theory. We have
preferred here, however, to use
``softer'' contour arguments. The understanding of the right
notion of contour for a specific model is in any case useful in many respects
(percolation, block dynamics,\ldots).

 However, unlike what
happens with sequential updating, not all these Gibbs measures need to
be stationary for the infinite volume PCA, the non-stationary ones
being periodic with period two. To conclude, we  exhibit a Gibbs
 measure which is not stationary for
the associated PCA.

\section{Shift invariant Probabilistic Cellular Automata}

Let $S$ be a finite set. For $\s \in S^{\Z^{d}}$, $\s = (\s_{i})_{i
\in \Z^{d}}$, and $\L \subset \Z^{d}$, we let $\s_{\L} \in S^{\L}$ its
restriction to $\L$. Sometimes, when no confusion arises, we omit the index $\L$
in $\s_{\L}$.

A time-homogeneous
Markov chain on $S^{\L}$ is determined, in law, by its transition
probabilities $P_{\L}(d\s|\eta)$. If $P_{\L}(d\s|\eta)$ is a product
measure, as a probability measure on $S^{\L}$, then we say that the
Markov chain is a Probabilistic Cellular Automaton. More explicitely
\[
P_{\L}(d\s|\eta) = \otimes_{i \in \L} P_{i}(d\s_{i}|\eta),
\]
and
\be{2.1}
P_{i}(\s_{i} = s| \eta) \equiv p_{i}(s|\eta), \ \ s \in S.
\ee
In the case $\L = \Z^{d}$, we omit the index $\L$ in
$P_{\L}(d\s|\eta)$. In this case, we say that a PCA is {\em shift
invariant} if, for every $i \in \Z^{d}$, $s \in S$, $\eta \in
S^{\Z^{d}}$, we have
\[
p_{i}(s|\eta) = p_{0}(s | \theta_{i}\eta),
\]
where $\theta_{i}$ is the shift in $\Z^{d}$: $(\theta_{i}
\eta)_{j} = \eta_{i+j}$ for every $j \in \Z^{d}$. A shift invariant
PCA is said to be {\em local} if, for each $s\in S$, the map $\eta
\ra p_{0}(s|\eta)$ is local, i.e. it depends on a finite number of
components of $\eta$.

From now on, all PCA's we consider in this paper satisfy the non
degeneration condition :
$$
p_{0}(s|\eta) >0, \quad \forall s \in S, \;\eta \in S^{\Z^{d}}.
$$
This means that we are
dealing with dynamics which can not contain a deterministic component.

In this paper we are mostly interested in stationary measures for
PCA's. For this purpose we recall the notion of Gibbs measure on
$S^{\Z^{d}}$. A shift invariant potential $\Phi$ is a family
$\{\Phi_{\L}: \L \subset \Z^{d}, \, |\L|<+\infty\}$ of maps
$\Phi_{\L}: S^{\L} \ra \R$ with the properties
\bi
\item[i.] For all $i \in \Z^{d}$, $\L \subset \Z^{d}$ finite:
\[
\Phi_{\L +i} = \Phi_{\L}\circ \theta_{i}.
\]
\item[ii.]
\[
\sum_{\L \ni 0} \|\Phi_{\L}\|_{\infty} < +\infty.
\]
\ei
Here and later $|\L|$ denotes the cardinality of $\L$.
Letting $H_{\L}(\s) = \sum_{A \cap \L \neq \emptyset}\Phi_{A}(\s)$ and
choosing $\tau \in S^{\Z^{d}}$,
also write $H_{\L}^{\tau}(\s_\L) = H_{\L}(\s_{\L}\tau_{\L^{c}})$, where
$\s_{\L}\tau_{\L^{c}}$ is the element of $S^{\Z^{d}}$ which coincides
with $\s$ on $\L$ and with $\tau$ on $\L^{c}$. The finite volume Gibbs
measure on $S^{\L}$ with boundary condition $\tau$ is given by
\[
\mu_{\L}^{\tau}(\s_\L) = \frac{\exp\left[-H_{\L}^{\tau}(\s_\L)
\right]}{Z_{\L}^{\tau}},
\]
where $Z_{\L}^{\tau}$ is the normalization factor. A probability
measure $\mu$ on $S^{\Z^{d}}$ is said to be Gibbsian for the
potential $\Phi$, and we write $\mu \in {\cal{G}}(\Phi)$ if for
every $\L \subset \Z^{d}$ finite and $\s \in S^{\Z^{d}}$
\[
\mu (\{\eta: \eta_{\L} = \s_{\L}\}| \eta_{\L^{c}} = \tau_{\L^{c}}) =
\mu_{\L}^{\tau}(\s_\L)
\]
for $\mu$-a.e. $\tau$. If $\mu$ is {\em shift-invariant}, i.e. $\mu
\circ \theta_{i} = \mu$ for all $i \in \Z^{d}$, then we write $\mu
\in {\cal{G}}_{s}(\Phi)$. More generally, we let ${\cal{P}}$
(resp. ${\cal{P}}_{s}$) be the set of probability measures (resp.
shift-invariant probability measures) on $S^{\Z^{d}}$.

Given $\L \subset \Z^{d}$, we denote by ${\cal{F}}_{\L}$ the
$\s$-field on $S^{\Z^{d}}$ generated by the projection $\s \ra \s_{\L}$.
For $\nu \in {\cal{P}}$, $\pi_{\L} \nu$ is the restriction of $\nu$
to ${\cal{F}}_{\L}$. We will use, for $\nu,\mu \in {\cal{P}}$, the
notion of local relative entropy:
\be{2.2}
h_{\L}(\nu|\mu) = \sum_{\s_{\L}} \pi_{\L}\nu(\s_{\L}) \log
\frac{\pi_{\L}\nu(\s_{\L})}{\pi_{\L}\mu(\s_{\L})}
\ee
with $\L \subset \Z^{d}$ finite, and of specific relative entropy
\be{2.3}
h(\nu|\mu) = \limsup_{\L \uparrow \Z^{d}} \frac{1}{|\L|}
h_{\L}(\nu|\mu)
\ee
where in the limit above $\L$ varies over hypercubes centered in the
origin. It is easily seen that $0 \leq h(\nu |\mu) \leq +\infty$. In
the case of $\mu \in {\cal{G}}_{s}(\Phi)$ for a potential $\Phi$,
in (\ref{2.2}) $\pi_{\L}\mu(\s_{\L})$ can be replaced by
$\mu_{\L}^{\tau}(\s_{\L})$, for an arbitrary $\tau$, without changing
the limit in (\ref{2.3}). Moreover, for $\mu \in {\cal{G}}_{s}(\Phi)$
and $\nu \in {\cal{P}}_{s}$, the limsup in (\ref{2.3}) is actually a
limit. In this case the Gibbs variational principle states that, for
$\nu \in {\cal{P}}_{s}$, $h(\nu|\mu) = 0$ if and only if $\nu \in
{\cal{G}}_{s}(\Phi)$; so $h(\nu|\mu)$ represents a notion of
(pseudo-) distance of $\nu$ from ${\cal{G}}_{s}(\Phi)$.

We now define a corresponding notion of specific relative entropy for
transition probabilities, that will be used to measure distance
between two dynamics. Let $P(d\s|\eta)$ and $Q(d\s|\eta)$ two
transition probabilities on $S^{\Z^{d}}$, and $\nu \in {\cal{P}}$. We
define
\[
{\cal H}_{\nu}(P|Q) = \limsup_{\L \uparrow \Z^{d}}\frac{1}{|\L|} \int
h_{\L}(P(\cdot|\eta)|Q(\cdot|\eta))\nu(d\eta).
\]
Clearly ${\cal H}_{\nu}(P|Q) \geq 0$. By conditioning to $\s$ the joint law
$Q_{\nu}(d\s,d\eta) \equiv P(d\s|\eta)\nu(d\eta)$
we obtain the backward transition probability,
that we denote by $\hat{P}_{\nu}(d\eta|\s)$. We also let $P\nu(d\s)$
be given by
\[
P\nu(A) = \int P(A|\eta) \nu(d\eta)
\]
for $A \subset S^{\Z^{d}}$ measurable. If $P\nu = \nu$ we say that
$\nu$ is stationary for $P(d\s|\eta)$.

Our first result concerns the entropy production for a PCA (cf. \cite{dp}). The
corresponding result in continuous time has appeared in \cite{ha}.
\begin{prpstn} \label{p2.1}
  Suppose $\mu $ is a stationary measure for
a shift invariant, local PCA with transition probability
$P(d\s|\eta)$. If  $\mu $ is also a shift invariant Gibbs measure
w.r.t. a certain potential $\Phi$ ( $\mu \in {\cal{G}}_{s}(\Phi)$),
then,
 for any shift invariant measure $\nu$,
\[
h(\nu|\mu) - h(P\nu |\mu) = {\cal H}_{\nu}(\hat{P}_{\nu}|\hat{P}_{\mu}).
\]
In particular, if $\nu \in {\cal{G}}_{s}(\Phi)$,
then $P\nu \in {\cal{G}}_{s}(\Phi)$, that is the set of
shift-invariant Gibbs measures w.r.t. the potential  $\Phi$ is stable
under the action of this PCA dynamics.
\end{prpstn}
\begin{proof}
Let $\L$ be a finite subset of $\Z^{d}$, and consider
\[
P_{\L}(\s|\eta) = \prod_{i \in \L}p_{i}(\s_{i}|\eta).
\]
This expression depends on the restriction of $\eta$ to a neighborhood
of $\L$, that we denote by $\overline{\L}$.

Consider now the measure $Q_{\nu}(d\s,d\eta)$ defined above. For $A,B
\subset \Z^{d}$ with $A$ finite, we denote by
$Q_{\nu}(\s_{A}|\eta_{B})$ the restriction to the $\s$-field generated
by the projection $(\s,\eta) \ra \s_{A}$ of the measure $Q$
conditioned to the $\s$-field generated by the projection $(\s,\eta)
\ra \eta_{B}$. So, e.g., $P_{\L}(\s|\eta) =
Q_{\nu}(\s_{\L}|\eta_{\Z^{d}}) \equiv Q_{\nu}(\s_{\L}|\eta)$,
independently of $\nu$. Similarly, $\hat{Q}_{\nu}(\eta_{A}|\s_{B})$ denotes
the time-reversed conditioning, so that
\be{2.4}
\pi_{\L} \hat{P}_{\nu}(\eta_{\L}|\s) = \hat{Q}_\nu(\eta_{\L}|\s).
\ee
For $C \subset \Z^{d}$ we will also use conditionings of the form
\[
\hat{Q}(\eta_{A}|\s_{B},\eta_{C}),
\]
with the obvious meaning.

A simple computation, using the fact that $P\mu = \mu$, yields
\begin{eqnarray*}
    h_{\overline{\L}}(\nu|\mu) - h_{\L}(P\nu |\mu) & = &  \\
    & = & \sum_{\s_{\L}} \pi_{\L}(P\nu)(\s_{\L})
    \sum_{\eta_{\overline{\L}}} \hat{Q}_{\nu}(\eta_{\overline{\L}} |
    \s_{\L}) \log \frac{\hat{Q}_{\nu}(\eta_{\overline{\L}} |
    \s_{\L})}{\hat{Q}_{\mu}(\eta_{\overline{\L}} |
    \s_{\L})} \\
    & = & E^{Q}\left[  \log \frac{\hat{Q}_{\nu}(\eta_{\overline{\L}} |
    \s_{\L})}{\hat{Q}_{\mu}(\eta_{\overline{\L}} |
    \s_{\L})} \right].
    \end{eqnarray*}
 Since
 \[
 h(\nu|\mu) - h(P\nu|\mu) = \lim_{\L \uparrow \Z^{d}}\frac{1}{|\L|} [
 h_{\overline{\L}}(\nu|\mu) - h_{\L}(P\nu |\mu)],
 \]
 then the conclusion follows provided we show (see (\ref{2.4}))
 \be{2.5}
  \lim_{\L \uparrow \Z^{d}}\frac{1}{|\L|} E^{Q}\left[  \log \frac{\hat{Q}_{\nu}(\eta_{\overline{\L}} |
    \s_{\L})}{\hat{Q}_{\nu}(\eta_{\overline{\L}} |
    \s)} \right] = 0
    \ee
    and
 \be{2.6}
 \lim_{\L \uparrow \Z^{d}}\frac{1}{|\L|} E^{Q}\left[
\log \frac{\hat{Q}_{\mu}(\eta_{\overline{\L}} |
    \s_{\L})}{\hat{Q}_{\mu}(\eta_{\overline{\L}} |
    \s)} \right] = 0.
    \ee
 Note that (\ref{2.6}) is a special case of (\ref{2.5}).

 Let now $\l^{\otimes}$ be the probability measure on $S^{\Z^{d}}$
 obtained by taking the infinite product of the uniform measure $\l$
 in $S$. We denote by $\l^{\otimes}(\s_{\L})$ the projection of
 $\l^{\otimes}$ on ${\cal{F}}_{\L}$. Let also
 $\{i_{1},\ldots,i_{|\overline{\L}|}\}$ be the lexicographic ordering
 of the elements of $\overline{\L}$; define $\overline{\L}_{k} =
 \{i_{1},\ldots,i_{k}\}$ for $1 \leq k \leq |\overline{\L}|$, and
 $\overline{\L}_{0} = \emptyset$. By the chain rule for conditional
 measures
 \be{2.7}
 \log \frac{\hat{Q}(\eta_{\overline{\L}}|
 \s_{\L})}{\l^{\otimes}(\eta_{\overline{\L}})} =
 \sum_{k=1}^{|\overline{\L}|}\log \frac{\hat{Q}(\eta_{i_{k}}|\s_{\L},
 \eta_{\overline{\L}_{k-1}})}{\l(\eta_{i_{k}})}.
 \ee
 Moreover, by shift invariance of $Q$
 \be{2.8}
 E^{Q} \left[ \log \frac{\hat{Q}(\eta_{i_{k}}|\s_{\L},
 \eta_{\overline{\L}_{k-1}})}{\l(\eta_{i_{k}})} \right] =
 E^{Q} \left[ \log \frac{\hat{Q}(\eta_{0}|\s_{\theta_{-i_{k}}\L},
 \eta_{\theta_{-i_{k}}\overline{\L}_{k-1}})}{\l(\eta_{0})} \right] .
 \ee
 Let $\Z^{d}_{-} = \{i \in \Z^{d}: i \prec 0\}$, where ``$\prec$'' is
 the lexicographic order. By the Shannon-Breiman-McMillan
 Theorem (\cite{ba}), for every $\e>0$ there are $A \subset \Z^{d}$, $B \subset \Z^{d}_{-}$
 finite such that if $A \subset V$ and $B \subset W \subset
 \Z^{d}_{-}$ then
 \be{2.9}
 \left| E^{Q}\left[ \log \frac{\hat{Q}(\eta_{0}| \s_{V},
 \eta_{W})}{\l(\eta_{0})} \right] -
  E^{Q}\left[ \log \frac{\hat{Q}(\eta_{0}| \s_A,
 \eta_B)}{\l(\eta_{0})} \right] \right| <\e.
 \ee
 Note that, if we take $\L$ large enough and $i_{k} \in \L$ is far
 enough from the boundary of $\L$, then $A\subset \theta_{-i_{k}}\L$,
 and $B \subset \theta_{-i_{k}}\overline{\L}_{k-1}$. For the other
 values of $i_{k} \in \overline{\L}$,
 \[
  E^{Q} \left[ \log \frac{\hat{Q}(\eta_{0}|\s_{\theta_{-i_{k}}\L},
 \eta_{\theta_{-i_{k}}\overline{\L}_{k-1}})}{\l(\eta_{0})} \right]
 \leq \log |S|,
 \]
 which is the upper bound for the entropy of any probability measure
 is $S$ with respect to $\l$. Summing all up
 \be{2.10}
 \lim_{\L \uparrow \Z^{d}}\frac{1}{|\L|} E^{Q}\left[\log \frac{\hat{Q}(\eta_{\overline{\L}}|
 \s_{\L})}{\l^{\otimes}(\eta_{\overline{\L}})}\right] =
  E^{Q}\left[ \log \frac{\hat{Q}(\eta_{0}| \s_{\Z^{d}},
 \eta_{\Z^d_-})}{\l(\eta_{0})} \right].
 \ee
 Exactly in the same way one shows that
 \be{2.11}
 \lim_{\L \uparrow \Z^{d}}\frac{1}{|\L|} E^{Q}\left[\log \frac{\hat{Q}(\eta_{\overline{\L}}|
 \s)}{\l^{\otimes}(\eta_{\overline{\L}})}\right] =
  E^{Q}\left[ \log \frac{\hat{Q}(\eta_{0}| \s_{\Z^{d}},
 \eta_{\Z^d_-})}{\l(\eta_{0})} \right].
 \ee
 Thus (\ref{2.10}) and (\ref{2.11}) establish (\ref{2.5}).
 \end{proof}

 Next result shows that the measures in ${\cal{P}}_{s}$ for which the
 entropy production is zero are exactly those in ${\cal{G}}_{s}(\Phi)$.
 This result goes back to \cite{ho}, where it has been proved for reversible
 systems in continuous time. The assumption of reversibility has been
 dropped in \cite{ku1}. In discrete-time, the proof for a special class of
 reversible PCA is given in \cite{kv}, Proposition 1. In the generality given here, the
 first proof was contained (but unpublished) in one of the authors'
 PhD Thesis (\cite{dp}). Later, a proof using general entropy
 arguments was given in \cite{mv}. In this
 paper we have preferred to emphasize the fact that the following
 result comes  from the precise entropy production
 formula presented in Proposition \ref{p2.1}.
\begin{prpstn} \label{p2.2}
 Under the same assumptions of Proposition \ref{p2.1}, suppose $\nu
 \in {\cal{P}}_{s}$ is such that
 \be{2.12}
 h(\nu|\mu) = h(P\nu|\mu)
 \ee
  (in particular, this happens when $\nu$ is stationary). Then $\nu
 \in {\cal{G}}_s(\Phi)$.
 \end{prpstn}
 \begin{proof}
 By what seen in Proposition \ref{p2.1}, (\ref{2.12}) amounts to
 \be{2.13}
 {\cal H}_\nu(\hat{P}_{\nu}|\hat{P}_{\mu}) = 0.
 \ee

 We now adapt a classical argument for Gibbs measures (see e.g. \cite{pr}, Th. 7.4).
 Let $V$ be a fixed hypercube and, for $k >0$,
 \[
 \partial_{k}V = \{ i \in V^{c} : \mbox{dist}(i,V) \leq k \},
 \]
 where dist($\cdot$) is the Euclidean distance.
 Take, now, a hypercube $\L_{m,k}$ that is obtained as disjoint union
 of $m^{d}$ translates of $V \cup \partial_{k}V$, say
 \[
 \L_{m,k} = \cup_{i=1}^{m^{d}} W_{i,k},
 \]
 where $W_{i,k} = T_{i}(V \cup \partial_{k}V)$, and $T_{i}$ is a
 suitable translation. We also write $V_{i} = T_{i}V$. Defining, for
 $i \in \{1,\ldots,m^{d}\}$
 \[
 B_{i,k} = W_{i,k} \setminus  V_{i}
 \]
 we have (we use the notations introduced in the proof of Proposition
 \ref{p2.1})
 \[
 \log
 \frac{\hat{Q}_{\nu}(\eta_{\L_{m,k}}|\s)}{\hat{Q}_{\mu}(\eta_{\L_{m,k}}|\s)}
 = \sum_{i=1}^{m^{d}} \log
 \frac{\hat{Q}_{\nu}(\eta_{V_{i}}|\eta_{B_{i,k}},\s)}
 {\hat{Q}_{\mu}(\eta_{V_{i}}|\eta_{B_{i,k}},\s)} + \log
 \frac{\hat{Q}_{\nu}(\eta_{B_{1,k}}|\s)}{\hat{Q}_{\mu}(\eta_{B_{1,k}}|\s)}.
 \]
 By positivity of relative entropy:
 \[
 E^{Q}\left[\log
 \frac{\hat{Q}_{\nu}(\eta_{B_{1,k}}|\s)}{\hat{Q}_{\mu}(\eta_{B_{1,k}}|\s)}
 \right] \geq 0
 \]
 so that
 \be{2.14}
 E^{Q}\left[\log
 \frac{\hat{Q}_{\nu}(\eta_{\L_{m,k}}|\s)}{\hat{Q}_{\mu}(\eta_{\L_{m,k}}|\s)}
 \right] \geq \sum_{i=1}^{m^{d}} E^{Q}\left[\log
 \frac{\hat{Q}_{\nu}(\eta_{V_{i}}|\eta_{B_{i,k}},\s)}
 {\hat{Q}_{\mu}(\eta_{V_{i}}|\eta_{B_{i,k}},\s)} \right].
 \ee
 By translation invariance of $Q$:
 \be{2.15}
 E^{Q}\left[\log
 \frac{\hat{Q}_{\nu}(\eta_{V_{i}}|\eta_{B_{i,k}},\s)}
 {\hat{Q}_{\mu}(\eta_{V_{i}}|\eta_{B_{i,k}},\s)} \right] =
 E^{Q}\left[\log
 \frac{\hat{Q}_{\nu}(\eta_{V}|\eta_{T_{i}^{-1}B_{i,k}},\s)}
 {\hat{Q}_{\mu}(\eta_{V}|\eta_{T_{i}^{-1}B_{i,k}},\s)} \right].
 \ee
 Moreover, since $B_{i,k} \uparrow V_{i}^{c}$ as $k \uparrow
 +\infty$, using again the Shannon-Breiman-McMillan Theorem, for
 each $\e>0$ we can choose $k$ large enough so that
 \be{2.16}
 \left| E^{Q}\left[\log
 \frac{\hat{Q}_{\nu}(\eta_{V}|\eta_{T_{i}^{-1}B_{i,k}},\s)}
 {\hat{Q}_{\mu}(\eta_{V}|\eta_{T_{i}^{-1}B_{i,k}},\s)} \right] -
 E^{Q}\left[\log
 \frac{\hat{Q}_{\nu}(\eta_{V}|\eta_{V^{c}},\s)}
 {\hat{Q}_{\mu}(\eta_{V}|\eta_{V^{c}},\s)}\right] \right| \leq \e.
 \ee
 Summing up (\ref{2.14}), (\ref{2.15}) and (\ref{2.16}), we get
 \be{2.17}
 \frac{1}{m^{d}} E^{Q} \left[ \log
 \frac{\hat{Q}_{\nu}(\eta_{\L_{m,k}}|\s)}{\hat{Q}_{\mu}(\eta_{\L_{m,k}}|\s)}
 \right] \geq E^{Q}\left[\log
 \frac{\hat{Q}_{\nu}(\eta_{V}|\eta_{V^{c}},\s)}
 {\hat{Q}_{\mu}(\eta_{V}|\eta_{V^{c}},\s)}\right] - \e.
 \ee
 But $m^{d}$ is proportional to $|\L_{m,k}|$, so, by (\ref{2.13})
 \[
 \lim_{m \ra +\infty} \frac{1}{m^{d}} E^{Q} \left[ \log
 \frac{\hat{Q}_{\nu}(\eta_{\L_{m,k}}|\s)}{\hat{Q}_{\mu}(\eta_{\L_{m,k}}|\s)}
 \right]  = 0.
 \]
 Thus, since $\e$ is arbitrary, (\ref{2.17}) yields
 \[
 E^{Q}\left[\log
 \frac{\hat{Q}_{\nu}(\eta_{V}|\eta_{V^{c}},\s)}
 {\hat{Q}_{\mu}(\eta_{V}|\eta_{V^{c}},\s)}\right] = 0
 \]
 that, by elementary properties of relative entropy, implies
 \be{2.18}
 \hat{Q}_{\nu}(\eta_{V}|\eta_{V^{c}},\s) =
 \hat{Q}_{\mu}(\eta_{V}|\eta_{V^{c}},\s) \ \ Q-\mbox{a.s.}
 \ee
 At this point we use Proposition 3.2 in \cite{ku2}, which implies that if
 (\ref{2.18}) holds for a $\mu \in {\cal{G}}_{s}(\Phi)$, then
$\nu(\eta_{V}|\eta_{V^{c}})= \mu(\eta_{V}|\eta_{V^{c}}) $ a.s. and
 then $\nu \in
 {\cal{G}}_{s}(\Phi)$ too. This completes the proof.
 \end{proof}

 \section{A class of reversible dynamics}

 In this section we introduce a class of reversible PCA's we will be dealing
 with in the rest of the paper, and give some general results on their
 stationary measures, resp. reversible measures. Let us remember that a
 PCA $P$ is called reversible if there exists at least one probability
 measure $\mu$ such that the Markov process with initial law $\mu$ and
 dynamics $P$ is reversible.

 We choose $S = \{-1,1\}$ as spin space. Consider a function
 $k:\Z^{d} \ra \R$ that is of finite range, i.e. there exists $R>0$
 such that $k(i) = 0$ for $|i|>R$, and symmetric, i.e. $k(i) = k(-i)$
 for every $i \in \Z^{d}$ (this last assumption being necessary to assure
 the reversibility of the PCA, cf \cite{kv}). Moreover, let $\tau \in \{-1,1\}^{\Z^{d}}$
 be a fixed configuration, that will play the role of boundary
 condition. For $\L \subset \Z^{d}$, we define the transition
 probability $P_{\L}^{\tau}(d\s |\eta) = \otimes_{i \in \L}
 P_{i}^{\tau}(d\s_{i}|\eta)$ by
 \be{PCA}
 P^{\tau}_{i}(\s_{i} = s|\eta) = p_{i}(s|\tilde{\eta}) = \frac{1}{2} \left[ 1+s
 \tanh (\b \sum_{i \in \Z^{d}} k(i-j)\tilde{\eta}_{j} + \beta h ),
 \right]
 \ee
 where $\tilde{\eta} = \eta_{\L} \tau_{\L^c}$ ;
 $h \in \R$, $\b >0$ are given parameters. According to \cite{kv},
 this particular form of $p_{i}$ is indeed the most general one for a
 shift invariant non degenerate local PCA on $\{-1,1\}^{\Z^{d}}$.

In the case $\L$ is a
 hypercube, we can also consider periodic boundary conditions. The
 associated transition probability is denoted by $P_{\L}^{per}$. In
 general, when $\L$ is finite, we write $P_{\L}^{\tau}(\s|\eta)$ in
 place of $P_{\L}^{\tau}(\{\s\}|\eta)$. In the case $\L = \Z^{d}$, the
 boundary condition $\tau$ plays no role, and will be omitted.

 In the rest of this section we establish some simple facts about
 stationary measures for these PCA's.
 \begin{prpstn} \label{p3.1}
 Let $\L \subset \Z^{d}$ finite, and $\tau \in \{-1,1\}^{\Z^{d}}$.
 Then the finite volume PCA with transition probability $P_{\L}^{\tau}(\s|\eta)$
 has a unique stationary measure $\nu_{\L}^{\tau}$ given by
 \[
 \nu_{\L}^{\tau}(\s) = \frac{1}{W_{\L}^{\tau}} \prod_{i \in \L} e^{\b
 h \s_{i}} \cosh \left[ \b \sum_{j \in \Z^{d}} k(i-j)
 \tilde{\s}_{j} + \b h \right] e^{\b \s_{i} \sum_{j \in \L^{c}}k(i-j)
 \tau_{j}},
 \]
 where, as before,
 $
 \tilde{\s}=\s_{\L} \tau_{\L^c}$,
 and $W_{\L}^{\tau}$ is the normalization. Moreover,
 $\nu_{\L}^{\tau}$ is reversible for $P_{\L}^{\tau}$.
 \end{prpstn}
\begin{proof}
It is clear that $P_{\L}^{\tau}(\s|\eta) >0$ $\forall \ \s,\eta$, so
 that the Markov chain with transition probability $P_{\L}^{\tau}$
 has a unique stationary measure. Thus, we only have to show that
 $\nu_{\L}^{\tau}$ is reversible, i.e.
 \be{3.1}
 P_{\L}^{\tau}(\s|\eta) \nu_{\L}^{\tau}(\eta) \equiv
 P_{\L}^{\tau}(\eta|\s) \nu_{\L}^{\tau}(\s).
 \ee
 Observe that, since $\s_i \in  \{-1,1\}, P_{\L}^{\tau}$ may be written in the form
 \[
 P_{\L}^{\tau}(\s|\eta) = \prod_{i \in \L} \frac{e^{\b \s_{i} \left(
 \sum_{j} k(i-j) \tilde{\eta}_{j} + h\right)}}{2 \cosh \left(\b\sum_{j}
 k(i-j) \tilde{\eta}_{j} +\b h\right)}.
 \]
 Thus (\ref{3.1}) amounts to
 \[
 \sum_{i \in \L} \sum_{j \in \Z^{d}} \s_{i} \tilde{\eta}_{j} k(i-j) +
 \sum_{i \in \L} \sum_{j \not\in \L} \eta_{i}\tau_{j}k(i-j) =
 \sum_{i \in \L} \sum_{j \in \Z^{d}} \eta_{i} \tilde{\s}_{j} k(i-j) +
 \sum_{i \in \L} \sum_{j \not\in \L} \s_{i}\tau_{j}k(i-j)
 \]
 which is easily checked.
 \end{proof}

 The above result on stationary measures for PCA's in finite volume,
 has an immediate consequence in infinite volume.

 \begin{prpstn} \label{p3.2}
 Let $\tau$ be any fixed boundary condition, and $\mu$ be any limit
 point of $\nu_{\L}^{\tau}$ as $\L \uparrow \Z^{d}$. Then $\mu$ is
 reversible for the infinite volume PCA defined in (\ref{PCA}), and $\mu$ is Gibbsian
 for the shift-invariant potential $\Phi$ given by
 \be{3.2}
 \begin{array}{lll}
     \Phi_{\{i\}}(\s_{i}) & = & - \b h \s_{i} \\
     \Phi_{U_{i}}(\s_{U_i}) & = & - \log \cosh \left[ \b  \sum_{j} k(i-j)
     \s_{j} + \b h \right] \\
     \Phi_{\L}(\s_{\L}) & = & 0 \ \ \text{otherwise},
 \end{array}
 \ee
 where $U_i =\{ j : k(i-j) \neq 0 \}$, that is finite by assumption.
 \end{prpstn}
 \begin{proof}
 Note that the finite volume Gibbs measure for $\Phi$ is
 \[
 \mu_{\L}^{\tau}(\s) = \frac{1}{Z_{\L}^{\tau}} \prod_{i: dist(i,\L) \leq R}
 \cosh \left[ \b  \sum_{j} k(i-j)
     \tilde{ \s}_{j} + \b h \right] e^{\b h \s_i},
 \]
 that differs from $\nu_{\L}^{\tau}$ only for boundary terms (and for
 the renormalization constant). The
 fact that the limit of $\nu_{\L}^{\tau}$ is Gibbsian for $\Phi$
 follows therefore from general facts on Gibbs measures (\cite{ge}). The
 reversibility of $\mu$ for the infinite volume PCA is obtained as
 follows. Let $f:\{-1,1\}^{\Z^{d}} \times \{-1,1\}^{\Z^{d}} \ra \R$ be
 a function which is local in both variables. For $\L$ large enough,
 reversibility of $\nu_{\L}^{\tau}$ yields
 \be{3.3}
 \sum_{\s,\tau}P_{\L}^{\tau}(\s|\eta) \nu_{\L}^{\tau}(\eta)
 f(\s,\eta) =
 \sum_{\s,\tau}P_{\L}^{\tau}(\eta|\s) \nu_{\L}^{\tau}(\s)
 f(\s,\eta) .
 \ee
 Note that, for $\L$ large enough, the boundary condition $\tau$ in
 $P_{\L}^{\tau}$ does not play any role in (\ref{3.3}). Thus,
 letting $\L \uparrow \Z^{d}$ in (\ref{3.3}) obtaining
 \be{3.4}
 \int P(d\s|\eta) \mu(d\eta) f(\s,\eta) = \int
 P(d\eta|\s)\mu(d\s)f(\s,\eta),
 \ee
 that establishes reversibility of $\mu$.
 \end{proof}
\begin{rm} \label{r3.1}
 Instead of fixed boundary conditions, one can choose periodic
 boundary
conditions. In this case, the finite volume measure defined by
 \[
 \nu_{\L}^{per}(\s) = \frac{1}{W_{\L}^{per}} \prod_{i \in \L} \cosh
 \left[ \b \sum_{j \in \Z^{d}} k(i-j) \tilde{\s}_{j} + \b h \right]
 e^{\b h \s_i}
 \]
 where $\tilde{\s}$ is the periodic continuation of $\s$, is the unique stationary
 reversible
 measure for $P_{\L}^{per}$. Remark that, in opposition to fixed
 boundary conditions, we now have  that $\nu_{\L}^{per}=
 \mu_{\L}^{per}$, which means that the finite volume stationary measure for the
 finite volume PCA is equal to the local specification of the
 associated Gibbs measure.
 \end{rm}

Moreover, the following result  gives a complete description of the
links between
the set of reversible measures for the PCA $P$ (which will be denoted
by $\cal{R}$), the set of stationary ones denoted by $\cal{S}$, the set
${\cal{G}}(\Phi)$ of Gibbs measures with respect to the potential
$\Phi$ defined by (\ref{3.2}), and their respective intersections with
the set of shift-invariant measures : ${\cal{R}}_{s}$,
${\cal{S}}_{s}$, ${\cal{G}}_{s}(\Phi)$.
 \begin{prpstn} \label{p3.3}
The reversible measures for the PCA $P$ defined in (\ref{PCA}) are
exactly those Gibbs measures w.r.t. $\Phi$ given in (\ref{3.2}) which
are also stationary :
\be{revstatgibbs}
{\cal{R}}={\cal{S}} \cap {\cal{G}}(\Phi).
\ee
Moreover, the subset of shift invariant reversible measures is equal
to the set of shift invariant stationary measures :
\be{revstat}
{\cal{R}}_{s} = {\cal{S}}_{s}.
\ee
 \end{prpstn}
\begin{proof}
The proof of the first assertion is based on the following proposition
 proved in \cite{kv} :

Let $P$ be a non degenerate local reversible PCA. Each reversible
measure $\mu$ for $P$ is Gibbs w.r.t. a certain potential $\Phi_P$.
 Reciprocally, any Gibbs measure w.r.t. $\Phi_P$ is either a reversible
 measure for $P$ or periodic of period two.

Since obviously ${\cal{R}} \subset {\cal{S}}$, the abovementioned proposition
implies ${\cal{R}} \subset {\cal{S}} \cap {\cal{G}}(\Phi)$. For the
reciprocal inclusion, since stationary measures can not be 2-periodic, a
stationary Gibbsian measure is necessarely a reversible one.

To prove the second assertion, note that by Proposition \ref{p3.2} and
 Remark \ref{r3.1},
 ${\cal S}_s \cap {\cal{G}}_{s}(\Phi) \ni \mu^{per}$. Thus
 Proposition  \ref{p2.2} applies, that is : ${\cal{S}}_{s}
\subset {\cal{G}}_{s}(\Phi)$.
On the other hand, from the first
 assertion: ${\cal R}_{s}={\cal{S}}_{s} \cap {\cal{G}}_{s}(\Phi)$.
 Then ${\cal{R}}_{s}={\cal{S}}_{s}$.
\end{proof}
\begin{remark} \label{r3.2}
The proof of Proposition \ref{p3.3} doesn't use the specific form of
 the PCA $P$. So equalities (\ref{revstatgibbs}) and (\ref{revstat})
hold as soon as Proposition
\ref{p2.2} and the abovementioned result of \cite{kv} apply, that is
 for the general class of local, shift invariant, non degenerate reversible
 PCA dynamics on $S^{\Z^d}$ for any $S$ finite.
 \end{remark}

 \section{Phase transition}

 In this section we show that for some reversible PCA it is indeed
 the case that not all Gibbs measures for the potential in (\ref{3.2}) are
 stationary. We treat those PCA defined in (\ref{PCA}) for which
 $k(i) = 0$ for $|i|>1$ (id est R=1), $h=0$ and $d=2$. Besides $\beta$, there are three
 parameters in the game: $k(0), k(e_{1})$ and $k(e_{2})$, where
 $e_{1},e_{2}$ are the basis vectors in $\R^{2}$. The first result
 concernes the existence of phase transition for the potential $\Phi$.

 \begin{prpstn} \label{p4.1}
 Assume $k(e_{1}) \neq 0$, $k(e_{2}) \neq 0$. Then there exists
 $\beta_{c} \in (0,+\infty)$ such that for $\beta > \beta_{c}$
 $|{\cal{G}}(\Phi)| >1$.
 \end{prpstn}
 \begin{proof}
 We divide the proof into different cases, depending on the signes
 of $k(0),k(e_{1}),k(e_{2})$. Note that the transformation
 $k(\cdot) \ra -k(\cdot)$ leaves invariant the potential $\Phi$.

 \noindent
 {\em Case 1: $k(0) \geq 0$, $k(e_{1})>0$, $k(e_{2})>0$.}

 For a given square $\L \subset \Z^{2}$, let $Cl_{m}(\L) = \{i \in
 \Z^{d} : dist(i,\L) \leq m\}$. Consider a fixed configuration $\s \in
 \{-1,+1\}^{\Z^{2}}$ such that $\s_{i} \equiv +1$ for $i \not\in
 \L$ ( $\s_{\L^c} \equiv +1$). Moreover let $\Z_{*}^{2} = \Z^{2} + (1/2,1/2)$. We recall the
 classical notion of {\em Peierls contour} associated to $\s$. We say that the segment joining
 two nearest neighbors $a,b \in \Z_{*}^{2}$ is marked if this segment
 separates two nearest neighbors $i,j \in \Z^{2}$ for which
 $\s_{i}\s_{j} =-1$. Marked segments form a finite family of closed, non
 self-intersecting,
 piecewise linear curves, that we call {Peierls contours}. Each
 segment of a contour $\g$ separates two nearest neighbors
 whose spins have different signes (they necessarily  belong
to  $Cl_{1}(\L)$). If $i,j$ are nearest neighbors
 separated by $\g$ and $\s_{i} = -1$ we write $i \in \partial^{-} \g$
 and $j \in \partial^{+} \g$. We call the union of the sets of sites $\partial^{-} \g$
and $\partial^{+} \g$ the {\em boundary of the contour } $\g$.
For each $i \in \L$ for which $\s_{i} = -1$,
 there is a {\em minimal Peierls contour} $\g$ { around} $i$, i.e. such
 that $i$ is in the interior of the closed curve $\g$.

 This notion of minimal contour is the one used for the Ising model.
 Here we have to modify it as follows. Two Peierls contours $\g,\g'$
 are called {\em adjacent} if their boundaries  have a common  point. We
 say that two Peierls contours $\g,\g'$ communicates if they belong to
 a sequence of Peierls contours $\g_{1},\ldots,\g_{n}$ such that for
 all $k$, $\g_{k}$ and $\g_{k+1}$ are adjacent. The relation of
 communicating is an equivalence relation. We call simply {\em
 contour} the union of the Peierls contours in an equivalence class.
 The minimal contour around $i$ with $\s_{i} = -1$ is the one formed
 by the equivalence class which contains the minimal Peierls contour
 around $i$. The boundary ($\partial^{+}$ or $\partial^{-}$) of a contour is
 simply the union of the boundaries of the Peierls contours that form
 it (see Fig. \ref{contour-dessin}).

\begin{figure} \label{contour-dessin}
\begin{center}
\includegraphics{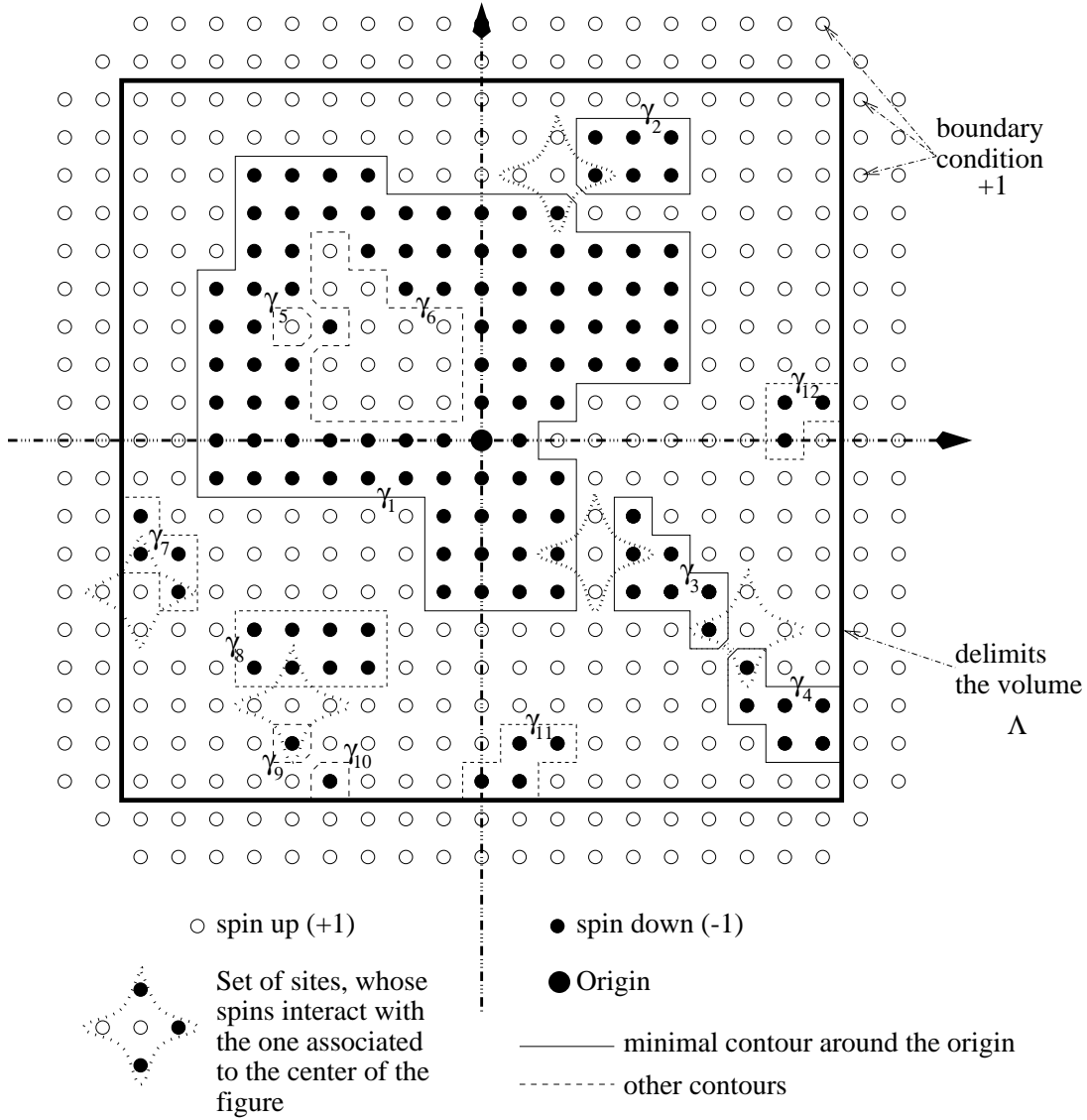}
\caption{Example of a configuration $\sigma$ on $Cl_{2}(\Lambda)$ such that $\sigma_0=-1$
  and \mbox{$\sigma_{\Lambda^c} \equiv +1$}. Drawing of its corresponding
  contours : %
$\gamma_1$ is the { \it minimal Peierls contour} around the origin ; %
\mbox{$( \gamma_1 \cup \gamma_2 \cup \gamma_3 \cup \gamma_4 )$} is the
{\it minimal contour} around the origin ({\it i.e.} the equivalence class of
$\gamma_1$) ; %
$\{ \gamma_5 , \gamma_6 \}, \{ \gamma_7 \}, \{ \gamma_8, \gamma_9,
\gamma_{10} \}, \{ \gamma_{11} \}, \{ \gamma_{12} \}$ are the other
equivalence classes.%
}
\end{center}
\end{figure}

 Let now $\mu_{\L}^+$ be the finite volume Gibbs measure with $+$
 boundary condition, that we write as follows:
\[
\mu^+_{\L}(\s) = \frac{1}{\Z_{\L}^+}\prod_{i \in Cl_1(\L)}
\frac{cosh(\b \sum_j k(i-j) \s_j^+)}{cosh(\b \sum_j k(i-j))} \quad
\text{ with } \s^+ = \s_\L(+1)_{\L^c}.
\]
We have modified the normalization for later convenience. A given $\s^+
\in \{-1,+1\}^{\Z^2}$ corresponds, as described above, to a collection of
contours $\Gamma = \{c_1,\ldots,c_m\}$. Each contour $c_i$ is a union
of Peierls contours. Peierls contours belonging to different $c_i$'s
do not communicate. We can write:
\[
\mu^+_{\L}(\s) = \frac{1}{\Z_{\L}^+} \prod_{k=1}^m F(c_k),
\]
where
\[
F(c_k) = \prod_{i \in \partial c_k}\frac{cosh(\b \sum_j k(i-j)
\s_j)}{cosh(\b \sum_j k(i-j))}
\]
and $\partial c_k = \partial^+ c_k \cup \partial^- c_k$. Observing
that if $\s_0 = -1$ then there is a contour around $0$, we have:
\[
\mu_{\L}^+ (\s_0 = -1) = \frac{1}{\Z_{\L}^+} \sum_{c_1 \text{
around } 0} F(c_1) \sum_{\Gamma \ni c_1} F(\Gamma \setminus c_1),
\]
where, for $\Gamma = c_1 \cup c_2 \cup \cdots \cup c_m$, we let
$F(\Gamma \setminus c_1) = \prod_{k=2}^m F(c_k)$. Note that, if
$\Gamma$ is a contour, $\Gamma \setminus c_1$ is also a contour,
that corresponds to the configuration obtained by flipping all the
spins $-1$ inside $c_1$ in the configuration associated to
$\Gamma$. It follows that
\[
 \sum_{\Gamma \ni c_1} F(\Gamma \setminus c_1) \leq Z_{\L}^+ \equiv
 \sum_{\Gamma} F(\Gamma),
 \]
 and therefore
 \be{4.1}
 \mu_{\L}^+ (\s_0 = -1) \leq \sum_{c_1 \text{
around } 0} F(c_1).
 \ee
 Now note that if $c_1$ is a contour and $i \in \partial c_1$,
 then the spins $\s_i, \s_{i \pm e_1}, \s_{i \pm e_2}$ do not have
 the same sign, so that
 \[
 \frac{cosh(\b \sum_j k(i-j)
 \s_j)}{cosh(\b \sum_j k(i-j))} \leq \frac{\cosh(\b A)}{\cosh(\b
 B)},
 \]
 where $B =\sum_j k(j)$, $A$ is the maximum value of $|\sum_j k(i-j)
 \s_j)|$ for $\s$ such that $\s_0, \s_{\pm e_1}, \s_{\pm e_2}$ do not have
 the same sign, and therefore $A<B$. Thus, we
 have to compare for a contour $c_1$, the cardinal of its boundary
$|\partial c_1|$ with its length denoted by $l(
 c_1)$. But remark that to any point of $\partial c_1$ correspond at
 most 4 marked segments on $c_1$. So, $l(c_1) \leq 4 |\partial c_1|$, and
 we have
 \[
 F(c_1) \leq \left[\frac{\cosh(\b A)}{\cosh(\b
 B)} \right]^{|\partial c_1|} \leq \left[\frac{\cosh(\b A)}{\cosh(\b
 B)} \right]^{l(c_1)/4}.
 \]
 On the other hand, for a given length $l$, it is easily checked
 that the number of contours around $0$ of length $l$ is
 bounded by
 $l^3 3^{l-1}$. Thus, by (\ref{4.1}),
 \[
  \mu_{\L}^+ (\s_0 = -1) \leq \sum_{l \geq 0} l^3 3^{l-1}
 \left[\frac{\cosh(\b A)}{\cosh(\b
 B)} \right]^{l/4}
 \]
 that goes to zero as $\b \uparrow +\infty$. Thus, taking $\b$
 large enough and letting $L \uparrow \Z^d$ in $\mu_{\L}^+$, we
 construct a Gibbs measure $\mu$ for which $\mu^+(\s_0 = -1) <1/2$.
 Simmetrically, taking minus boundary conditions, we obtain a
 Gibbs measure $\mu^-$ for which $\mu^-(\s_0 = -1) >1/2$, and this
 proves phase transition.

\noindent
{\em Case 2: $k(0) < 0$, $k(e_1) >0$, $k(e_2)>0$}.

Define
\[
k^*(i) = \left\{
\begin{array}{ll}
k(i) & \text{for } i \neq 0 \\
-k(0) & \text{for } i=0,
\end{array}
\right.
\]
and let $\Phi^*$ be the associated potential. Consider also the map
$T:\{-1,1\}^{\Z^2} \ra \{-1,1\}^{\Z^2}$ given by
\[
(T\s)_i = \left\{
\begin{array}{ll}
\s_i & \text{for } i \in \Z^2_e \\
-\s_i & \text{for } i \in \Z^2_o.
\end{array}
\right.
\]
To stress dependence on the potential $\Phi$ we write
$\mu_{\L,\Phi}^{\tau}$ for $\mu_{\L}^{\tau}$. It is easily seen that
\[
\mu_{\L,\Phi}^{\tau}(\s) = \mu_{\L,\Phi^*}^{T\tau}(T\s),
\]
so that the map $\mu \ra \mu \circ T$ is a bijection between
$\G(\Phi)$ and $\G(\Phi^*)$. The conclusion follows from the fact that
$|\G(\Phi^*)|>1$, as seen in case 1.

\noindent
{\em Case 3: $k(0) \geq 0$, $k(e_1)>0$, $k(e_2)<0$}.

This case is treated as case 2, with the following choices:
\[
k^*(i) = \left\{
\begin{array}{ll}
k(i) & \text{for } i \neq e_2 \\
-k(e_2) & \text{for } i=e_2,
\end{array}
\right.
\]
and
\[
(T\s)_i = \left\{
\begin{array}{ll}
\s_i & \text{for } i = (x,y) \text{ with } y \text{ even}\\
-\s_i & \text{otherwise}.
\end{array}
\right.
\]
the proof is now completed.
 \end{proof}
\begin{remark} \label{r4.1}

The special case $k(0)=0$  was already treated in  \cite{kv} example
2 (for $k(e_1)=k(e_2)=1$) , where a remarkable relation with Ising
model was pointed out. We recall here in some more generality the
principal steps of the argumentation :

let $\Z^2_o = \{ (x,y) \in \Z^2 : x+y \text{ is odd}\}$, $\Z^2_e =
\Z^2 \setminus \Z^2_o$ and, similarly, $\L_o = \L \cap \Z^2_o$,
$\L_e = \L \cap \Z^2_e$. Note that since $k(0)=0, \s_{\L_o}$ and $\s_{\L_e}$ are
independent under $\mu_{\L}^{\tau}$, i.e. $\mu_{\L}^{\tau} =
\mu_{\L_e}^{\tau} \otimes \mu_{\L_o}^{\tau}$. Consider the
following anisotropic Ising model on $\{-1,1\}^{\L}$:
\[
\rho_{\L}^{\tau}(\s) = \frac{1}{N_{\L}^{\tau}} \exp \left[ \b \sum_{i
\in \L} \left( k(e_1) \s_i \tilde{\s}_{i+e_1} + k(e_2) \s_i
\tilde{\s}_{i+e_2}\right) \right],
\]
where $N_{\L}^{\tau}$ is the normalization and $\tilde{\s} = \s_{\L}
\tau_{\L^c}$. Restricting this measure to the sites in $\L_e$ we
obtain
\begin{eqnarray*}
\pi_{\L_e} \rho_{\L}^{\tau}(\s_{\L_e}) & = & \sum_{\s_{\L_o}}
\rho_{\L}^{\tau}(\s) \\
 & = & \frac{2}{N_{\L}^{\tau}} \prod_{i \in \L_o} \cosh \left[ \b
\sum_j k(i-j) \tilde{\s}_j \right] \\
 & = & \mu_{\L_o}^{\tau}(\s_{\L_e}).
\end{eqnarray*}
Therefore, phase transition for $\G(\Phi)$ follows from phase
transition for the Ising model:\\
  since $\rho^-(\s_0=+1) < \frac{1}{2}<
\rho^+(\s_0=+1)$, the restrictions $\pi_{\L_e} \rho^-$ and $\pi_{\L_e}
\rho^+$ are different, and then
$$
(\mu^+_{\Z^2_o} = \pi_{\L_e} \rho^+ ) \neq  (\pi_{\L_e} \rho^- =
\mu^-_{\Z^2_o}).
$$
\end{remark}
We now show that, in certain cases, there are elements in $\G(\Phi)$
that are not stationary.

\begin{prpstn} \label{p4.2}
Suppose $k(0) \leq 0$, $k(e_1)<0$, $k(e_2)<0$, and let $\mu^+$ be the
Gibbs mesure corresponding to plus boundary conditions. Suppose $\b$
is large enough so that $\mu^+ \neq \mu^-$. Then $\mu^+$
is not stationary.
\end{prpstn}
\begin{proof}
We first observe that the transformation $k(\cdot) \ra -k(\cdot)$ do
not change the elements of $\G(\Phi)$, but it does change the
dynamics. We recall few basic notions on stochastic ordering. Given
$\s,\eta \in \{-1,1\}^{\Z^2}$, we say that $\s \leq \eta$ if $\s_i
\leq \eta_i$ for every $i \in \Z^2$. Monotonicity of functions
$\{-1,1\}^{\Z^2} \ra \R$ is defined with respect to this partial
order. Finally, for $\nu,\mu$ probabilities on $\{-1,1\}^{\Z^2}$, we
say that $\nu \leq \mu$ if $\int f d\nu \leq \int f d\mu$ for every
increasing $f$.

The key observation consists in the fact that, under our assumptions
on $k(\cdot)$, the transition probability $P(d\s|\eta)$ is decreasing,
i.e.
\[
\nu \leq \mu \text{ implies } P\nu \geq P\mu.
\]
This follows from the facts that $p_0(1|\eta)$ is decreasing in
$\eta$, while $p_0(-1|\eta)$ is increasing in $\eta$ (see \cite{lms}
or \cite{ls} for
details). Let now $\mu^0$ be a limit point of the sequence
$\nu_{\L}^{per}$ defined in Remark \ref{r3.1}. By using the criterion in
\cite{l}, Th. II 2.9,
it is easy to check that $\nu_{\L}^{per} \leq \nu_{\L}^+$ for every
$\L$, and so $\mu^0 \leq \mu^+$. Moreover, $ 0= \mu^0 (\s_0=-1) <
\frac{1}{2} <\mu^+ (\s_0=-1)$. So $\mu^0 \neq \mu^+$. On the other hand,
 by Proposition \ref{p3.2},
$\mu^0$ is stationary. Therefore $P\mu^+ \leq P\mu^0 = \mu^0 < \mu^+$,
which completes the proof.
\end{proof}


\bigskip


\begin{acknowledgement} : P.-Y. Louis thanks, for their kind hospitality,
the Mathematics'
Departement of Padova University and the Interacting Random Systems
group of Weierstrass Institute for Applied Analysis and Stochastics in Berlin, where part
of this work was done.
\end{acknowledgement}


\end{document}